\documentclass[12pt]{article}
\usepackage{latexsym,amsfonts,amssymb}
\usepackage[dvips]{graphicx}
\usepackage{color}
\usepackage{colordvi,multicol}

\usepackage{mathrsfs}
\usepackage{amsmath}
\usepackage{amssymb,amsmath,amsfonts,mathrsfs}

\usepackage{latexsym,amsfonts,amssymb}
\usepackage[dvips]{graphicx}
\usepackage{amscd}
\usepackage{graphicx}
\usepackage{subfigure}
\usepackage[mathlines]{lineno}

\newtheorem{thm}{Theorem}[section]

\newtheorem{pro}[thm]{Proposition}

\date{}
\begin{document}

\title{\bf A study on the Poisson, geometric and Pascal distributions motivated by Chv\'{a}tal's conjecture}
 \author{Fu-Bo Li$^{1}$, Kun Xu$^{2}$, Ze-Chun Hu$^{1,}$\footnote{Corresponding author: zchu@scu.edu.cn}\\ \\
  {\small $^{1}$ College of Mathematics, Sichuan  University,
 Chengdu 610065, China}\\
{\small $^{2}$ Deyang Foreign Languages School,
 Deyang 618099, China}}
\maketitle

\begin{abstract}

  Let $B(n,p)$ denote a binomial random variable with parameters $n$ and $p$. Vas\v{e}k Chv\'{a}tal conjectured that for any fixed $n\geq 2$, as $m$ ranges over $\{0,\ldots,n\}$, the probability $q_m:=P(B(n,m/n)\leq m)$ is the smallest when  $m$ is closest to $\frac{2n}{3}$. This conjecture has been solved recently.  Motivated by this conjecture, in this paper, we consider the corresponding minimum value problem on the probability that a random variable is not more than its expectation, when its distribution is the Poisson distribution, the geometric distribution or the Pascal distribution.

\end{abstract}



\section{Introduction}

Let $B(n,p)$ denote a binomial random variable with parameters $n$ and $p$. Janson (2021) introduced the following conjecture suggested by Va\v{s}k Chv\'{a}tal.

{\bf Conjecture 1} (Chv\'{a}tal). For any fixed $n\geq 2$, as $m$ ranges over $\{0,\ldots,n\}$, the probability $q_m:=P(B(n,m/n)\leq m)$ is the smallest when
 $m$ is closest to $\frac{2n}{3}$.

 As to the probability of a binomial random variable exceeding its expectation, we refer to Doerr (2018), Greenberg  and Mohri (2014), and Pelekis and Ramon (2016).   Janson (2021) proved that  Conjecture 1 holds for large $n$. Barabesi et al. (2023) and Sun (2021) gave an affirmative answer to Conjecture 1  for general $n\geq 2$.

Motivated by Conjecture 1,  we will consider the corresponding minimum value problem on the probability that a random variable is not more than its expectation, when its distribution is the Poisson distribution, the geometric distribution or the Pascal distribution in Sections 2, 3, 4, respectively.

Sun et al. \cite{SHS23} considered the minimum value problem for the Gamma distribution among other things.


\section{The Poisson distribution}\setcounter{equation}{0}

Let $X_{\lambda}$ denote a Poisson random variable  with parameter $\lambda\,(\lambda>0)$. Then
\begin{eqnarray*}
P(X_{\lambda}\leq E[X_{\lambda}])=P(X_{\lambda}\leq \lambda)=\sum_{k=0}^{[\lambda]}\frac{\lambda^ke^{-\lambda}}{k!}=
\left(1+\lambda+\frac{\lambda^2}{2!}+\cdots+\frac{\lambda^{[\lambda]}}{[\lambda]!}\right)
e^{-\lambda}.
\end{eqnarray*}
Hereafter, for a real number $a, [a]$ stands for the biggest integer which is not more than $a$.

The main result of this section is

\begin{pro}\label{pro-2.1}
The function $P(X_{\lambda}\leq E[X_{\lambda}])$ has no minimum value on $(0,\infty)$, but \begin{eqnarray}\label{pro-2.1-a}
\inf_{\lambda\in (0,\infty)}P(X_{\lambda}\leq \lambda)=\lim_{\lambda\uparrow 1}P(X_{\lambda}\leq \lambda)=\frac{1}{e}.
\end{eqnarray}
\end{pro}

\noindent {\bf Proof.} For $\lambda\in (0,1)$, we have
$$
P(X_{\lambda}\leq \lambda)=e^{-\lambda}.
$$
It follows that
\begin{eqnarray}\label{pro-2.1-b}
\inf_{\lambda\in (0,1)}P(X_{\lambda}\leq \lambda)=\lim_{\lambda\uparrow 1}e^{-\lambda}=\frac{1}{e}.
\end{eqnarray}

In the first version of this paper (see \cite{XLH22}), we used pure analysis method to give a long proof of
$$
\inf_{\lambda\in [1,\infty)}P(X_{\lambda}\leq \lambda)>\frac{1}{e}.
$$
The following is a much shorter alternagive proof to (\ref{pro-2.1-a}) suggested by the referee to the first version.

By the central limit theorem,
$$
\lim_{\lambda\to\infty}P(X_{\lambda}\leq \lambda)=\lim_{\lambda\to\infty} P(\lambda^{-1/2}(X_{\lambda}- \lambda)\leq 0)=\frac{1}{2}.
$$
Therefore, by virtue of (\ref{pro-2.1-b}),  there is an integer $j$ such that
$$
\inf_{\lambda>j}P(X_{\lambda}\leq \lambda)>\frac{1}{e}=\inf_{\lambda\in (0,1)}P(X_{\lambda}\leq \lambda).
$$
Since $P(X_{\lambda}\leq x)$ is decreasing in $\lambda$ for fixed $x$, it follows that
$$
\inf_{\lambda>0}P(X_{\lambda}\leq \lambda)=\inf_{0<\lambda\leq j}P(X_{\lambda}\leq \lambda)=\inf_{0<\lambda\leq j}P(X_{\lambda}\leq [\lambda])=\min_{k=0,1,\ldots,j-1}P(X_{1+k}\leq k).
$$
Next, since $n-B(n,(k+1)/n)$ and $B(n,1-(k+1)/n)$ have the same distribution, one obtains
\begin{eqnarray*}
P(X_{k+1}\geq k+1)&=&\lim_{n\to\infty}P(B(n,(k+1)/n)\geq E[B(n,(k+1)/n)])\\
&=&\lim_{n\to\infty}P(B(n,1-(k+1)/n)\leq E[B(n,1-(k+1)/n]).
\end{eqnarray*}
In addition, by \cite[Theorem 2]{BPR23}, for large $n$, $P(B(n,1-(k+1)/n)\leq E[B(n,1-(k+1)/n])$ is decreasing in $k\in \{0,1,\ldots,j-1\}$. Taking the limit, $P(X_{k+1}\geq k+1)$ is decreasing in $k$, so that
$P(X_{1+k}\leq k)=1-P(X_{1+k}\geq k+1)$ is increasing in $k$. Hence
$$
\min_{k=0,1,\ldots,j-1}P(X_{1+k}\leq k)=P(X_1\leq 0)=\frac{1}{e}.
$$
The proof is complete.\hfill\fbox

\section{The geometric distribution}\setcounter{equation}{0}

Let $Y_p$ denote a geometric random variable with parameter $p \,(0<p\leq1)$.  Then we have
\begin{eqnarray*}
P(Y_p\leq E[Y_p])&=&P(Y\leq 1/p)\\
&=&\sum_{k=1}^{[1/p]}p(1-p)^{k-1}\nonumber\\
&=&\left(p+p(1-p)+\cdots+p(1-p)^{[1/p]-1}\right)\nonumber\\
&=&1-(1-p)^{[1/p]}.
\end{eqnarray*}

Define
\begin{eqnarray}\label{sec-3-a}
f(p):=1-(1-p)^{[1/p]},\quad 0<p\leq1.
\end{eqnarray}
The main result of this section is

\begin{pro}\label{pro-3.1}
The function $f(p)$ has no minimum value on $(0,1]$, but
\begin{eqnarray}\label{pro-3.1-a}
\inf_{p\in (0,1]}f(p)=\lim_{p\downarrow \frac{1}{2}}f(p)=\frac{1}{2}.
\end{eqnarray}

\end{pro}

\noindent {\bf Proof.} Let $x$ be a positive integer. For any $1/p\in[x,x+1)$, we have
$$
f(p)=1-(1-p)^x.
$$
Then for $p\in(\frac{1}{x+1},\frac{1}{x})$,
$$
f'(p)=x(1-p)^{x-1}>0,
$$
which implies that the function $f(p)$ is strictly increasing on the interval $(\frac{1}{x+1},\frac{1}{x}]$. Thus  we have
\begin{eqnarray*}
\inf_{p\in (\frac{1}{x+1},\frac{1}{x}]}f(p)
=\lim_{p\downarrow \frac{1}{x+1}}f(p)=1-\left(1-\frac{1}{x+1}\right)^x.
\end{eqnarray*}

Define a sequence $\{a_n\}$ as follows:
\begin{eqnarray}
a_n:=1-\left(1-\frac{1}{n+1}\right)^n,\ \forall n\geq 1.
\end{eqnarray}

Since the sequence $(1+1/n)^n$ is increasing and
$$
a_n=1-\left(1-\frac{1}{n+1}\right)^n=1-\frac{1}{(1+1/n)^n},
$$
it follows that $\min_{n\geq 1}a_n=a_1=1/2$.

By the analysis above, we know that $f(p)$ has no minimum value on $(0,1]$, but its infimum satisfies
\begin{eqnarray*}
\inf_{p\in (0,1]}f(p)=\lim_{p\downarrow \frac{1}{2}}f(p)=a_1=\frac{1}{2}.
\end{eqnarray*}
Hence (\ref{pro-3.1-a}) holds. The proof is complete.\hfill\fbox

\section{The Pascal distribution}\setcounter{equation}{0}

Let $B^*(r,p)$ denote a Pascal random variable with parameters $r$ ($r$ is a positive integer) and $p\, (0<p\leq1)$ such that
\begin{eqnarray*}
P(B^*(r,p)=j)=
\left(
\begin{array}{l}
j-1\\
r-1
\end{array}
\right)(1-p)^{j-r}p^r,\ j=r,r+1,\cdots.
\end{eqnarray*}
We know that the expectation  is $E[B^*(r,p)]=r/p$.

An anonymous referee to the first version of this paper posed the following {\bf conjecture}:
\begin{eqnarray}\label{4.1}
\inf_{0<p\leq 1}P\left(B^*(r,p)\leq r/p\right)&=&\inf_{n\geq 1}P\left(B\left(n+r-1,\frac{r}{n+r}\right)\geq r\right)\nonumber\\
&=&\inf_{n\geq r}P\left(B\left(n,\frac{r}{n+1}\right)\geq r\right)\nonumber\\
&=&P\left(B\left(r,\frac{r}{r+1}\right)\geq r\right)\nonumber\\
&=&\left(\frac{r}{r+1}\right)^r,
\end{eqnarray}
and presented the proof to the first equality in (\ref{4.1}) as follows:
\begin{eqnarray*}
P\left(B^*(r,p)\leq r/p\right)=P\left(B^*(r,p)\leq [r/p]\right)=P\left(B([r/p],p)\geq r\right).
\end{eqnarray*}
For each $p\in (\frac{r}{n+r},\frac{r}{n+r-1}]$, one obtains
\begin{eqnarray*}
P(B([r/p],p)\geq r)&=&P(B(n+r-1,p)\geq r)\\
&\geq &P\left(B\left(n+r-1,\frac{r}{n+r}\right)\geq r\right)\\
&=&\inf_{q\in (\frac{r}{n+r},\frac{r}{n+r-1}]}P(B([r/q],q)\geq r).
\end{eqnarray*}
It follows that the first equality in (\ref{4.1}) holds.


Define a function
\begin{eqnarray*}\label{4.4}
f_r(p):=\sum_{k=r}^{[r/p]}\dbinom{k-1}{r-1}p^r(1-p)^{k-r},\ p\in (0,1].
\end{eqnarray*}
Then Conjecture (\ref{4.1}) becomes
\begin{eqnarray}\label{4.4}
\inf_{0<p\leq 1}f_r(p)=\left(\frac{r}{r+1}\right)^r.
\end{eqnarray}

When $\frac{r}{p}\in [n,n+1) (n\geq r, n\in \mathbb{N}$), we have
\begin{eqnarray*}
f_r(p):=\sum_{k=r}^{n}\dbinom{k-1}{r-1}p^r(1-p)^{k-r},\ p\in \left(\frac{r}{n+1},\frac{r}{n}\right].
\end{eqnarray*}
For $k=r,\ldots,n$, define a function
$$
g_k(p):=\dbinom{k-1}{r-1}p^r(1-p)^{k-r},\ p\in \left(\frac{r}{n+1},\frac{r}{n}\right].
$$
By the condition that $\frac{r}{k}\geq \frac{r}{n}$, we get
$$
\frac{dg_k(p)}{dp}=\dbinom{k-1}{r-1}
\frac{kp^r(1-p)^{k-r}(\frac{r}{k}-p)}{p(1-p)}>0,\ \forall p\in \left(\frac{r}{n+1},\frac{r}{n}\right).
$$
It follows that for any $k=r,\ldots,n$, the function $g_k(p)$ is strictly increasing on $\left(\frac{r}{n+1},\frac{r}{n}\right]$. Then
\begin{eqnarray*}
\inf_{p\in (\frac{r}{n+1},\frac{r}{n}]}g_k(p)
=\lim_{p\downarrow \frac{r}{n+1}}g_k(p)=\dbinom{k-1}{r-1}\left(\frac{r}{n+1}\right)^r
\left(1-\frac{r}{n+1}\right)^{k-r}.
\end{eqnarray*}
Thus
\begin{eqnarray*}
\inf_{p\in (\frac{r}{n+1},\frac{r}{n}]}f_r(p)
=\lim_{p\downarrow \frac{r}{n+1}}f_r(p)=\sum_{k=r}^{n}\dbinom{k-1}{r-1}
\left(\frac{r}{n+1}\right)^r\left(1-\frac{r}{n+1}\right)^{k-r}.
\end{eqnarray*}

Define a number sequence
\begin{eqnarray}\label{4.5}
a_r(n)&:=&\sum_{k=r}^{n}\dbinom{k-1}{r-1}\left(\frac{r}{n+1}\right)^r
\left(1-\frac{r}{n+1}\right)^{k-r}\nonumber\\
&=&\left(\frac{r}{n+1}\right)^r\sum_{k=r}^{n}\dbinom{k-1}{r-1}
\left(1-\frac{r}{n+1}\right)^{k-r},\quad \forall n\geq r.
\end{eqnarray}
Then (\ref{4.4}) (or equivalently, Conjecture (\ref{4.1})) becomes
\begin{eqnarray}\label{4.6}
\min_{n\geq r}a_r(n)=a_r(r).
\end{eqnarray}

Up to now, we have not found the unified method to prove (\ref{4.6}) for any $r\geq 2$.
In the following two subsections, we will give the proofs for the cases $r=2,3$, respectively.

\subsection{Case $r=2$}

In this case, we have
\begin{eqnarray}\label{r=2-a}
a_2(n)&=&\left(\frac{2}{n+1}\right)^2\sum_{k=2}^{n}(k-1)\left(1-\frac{2}{n+1}\right)^{k-2}\nonumber\\
&=&\left\{
\begin{array}{cl}
\frac{4}{9},&\mbox{if}\ n=2,\\
1-\frac{3n-1}{n+1}\left(\frac{n-1}{n+1}\right)^{n-1}, & \mbox{if}\ n\geq 3.
\end{array}
\right.
\end{eqnarray}

\begin{pro}\label{pro-4.1}
\begin{eqnarray}\label{pro-4.1-a}
\min_{n\geq 2}a_2(n)=a_2(2).
\end{eqnarray}
\end{pro}
{\bf Proof.} Define a number sequence
\begin{eqnarray}\label{r=2-b}
b_2(n):=\frac{3n-1}{n+1}\left(\frac{n-1}{n+1}\right)^{n-1},\ n\geq 3,
\end{eqnarray}
and a function
\begin{eqnarray}\label{r=2-c}
g_2(x):=\frac{3x-1}{x+1}\left(\frac{x-1}{x+1}\right)^{x-1},\ x\geq 3.
\end{eqnarray}
Then we have
\begin{eqnarray}\label{r=2-d}
g_2'(x)=\frac{3x-1}{x+1}\left(\frac{x-1}{x+1}\right)^{x-1}\left(\frac{3}{3x-1}+\frac{1}{x+1}+\ln \frac{x-1}{x+1}\right) .
\end{eqnarray}

Define a function
\begin{eqnarray}\label{r=2-e}
h_2(x):=\frac{3}{3x-1}+\frac{1}{x+1}+\ln \frac{x-1}{x+1},\ x\geq 3.
\end{eqnarray}
Then we have
\begin{eqnarray*}
h_2'(x)=\frac{4(3x^2-2x+3)}{(x+1)^2(x-1)(3x-1)^2}.
\end{eqnarray*}
It is easy to check that $3x^2-2x+3>0, \forall x\in (-\infty,\infty)$. Then we get
$$
h_2'(x)>0,\quad \forall x\geq 3,
$$
which implies that $h_2(x)$ is strictly increasing on $[3,\infty)$. Since
\begin{eqnarray*}
\ \lim_{x\rightarrow +\infty}h_2(x)
=\lim_{x\rightarrow +\infty}\left(\frac{3}{3x-1}+\frac{1}{x+1}+\ln \frac{x-1}{x+1}\right)=0,
\end{eqnarray*}
 we get that $h_2(x)<0,\forall x\geq 3$. By (\ref{r=2-d}) and (\ref{r=2-e}), we know that $g_2'(x)<0,\forall x\geq 3$. Thus, by (\ref{r=2-b}) and (\ref{r=2-c}), we get
$$
\max_{n\geq 3}b_2(n)=b_2(3)=\frac{1}{2}<\frac{5}{9}.
$$
By (\ref{r=2-a}) and (\ref{r=2-b}), we get
$$
\min_{n\geq 2}a_2(n)=a_2(2)=\frac{4}{9},
$$
i.e. (\ref{pro-4.1-a}) holds. The proof is complete.\hfill\fbox

\subsection{Case $r=3$}

In this case, we have
\begin{eqnarray}\label{r=3-a}
a_3(n)&=&\left(\frac{3}{n+1}\right)^3\sum_{k=3}^{n}\dbinom{k-1}{2}
\left(1-\frac{3}{n+1}\right)^{k-3}\nonumber\\
&=&\left\{
\begin{array}{cl}
\frac{27}{64},&\mbox{if}\ n=3,\\
1-\frac{17n^2-29n+8}{2(n+1)^2}\left(\frac{n-2}{n+1}\right)^{n-2}, & \mbox{if}\ n\geq 4.\ \
\end{array}
\right.
\end{eqnarray}

\begin{pro}\label{pro-4.2}
\begin{eqnarray}\label{pro-4.2-a}
\min_{n\geq 3}a_3(n)=a_3(3).
\end{eqnarray}
\end{pro}
{\bf Proof.} Define a number sequence
\begin{eqnarray}\label{r=3-b}
b_3(n):=\frac{17n^2-29n+8}{2(n+1)^2}\left(\frac{n-2}{n+1}\right)^{n-2},
\ n\geq 4.
\end{eqnarray}
We will prove
\begin{eqnarray}\label{r=3-c}
\max_{n\geq 4}b_3(n)=b_3(4),
\end{eqnarray}
which implies that
$$
\min_{n\geq 4}a_3(n)=a_3(4)=1-\frac{328}{390625}>\frac{27}{64}.
$$
and thus
$$
\min_{n\geq 3}a_3(n)=a_3(3),
$$
i.e. (\ref{pro-4.2-a}) holds.

In the following, we will prove (\ref{r=3-c}).  Define a function
\begin{eqnarray}\label{r=3-d}
g_3(x):=\frac{17x^2-29x+8}{2(x+1)^2}\left(\frac{x-2}{x+1}\right)^{x-2},
\ x\geq 4.
\end{eqnarray}

We have
\begin{eqnarray}\label{r=3-e}
g_3'(x)=\frac{17x^2-29x+8}{2(x+1)^2}\left(\frac{x-2}{x+1}\right)^{x-2}
\left[\frac{34x-29}{17x^2-29x+8}+\frac{1}{x+1}+\ln\frac{x-2}{x+1}\right],
\ x\geq 4.
\end{eqnarray}

Define a function
\begin{eqnarray}\label{r=3-f}
h_3(x):=\frac{34x-29}{17x^2-29x+8}+\frac{1}{x+1}+\ln\frac{x-2}{x+1},
\ x\geq 4.
\end{eqnarray}
Then for any $x\geq 4$,
\begin{eqnarray*}
h_3'(x)
=\frac{27(17x^4-57x^3+105x^2-91x+54)}{(x+1)^2(x-2)(17x^2-29x+8)^2}>0.
\end{eqnarray*}
Hence $h_3(x)$ is strictly increasing on $[4,\infty)$. Since
$
\lim_{x\to\infty}h_3(x)=0,
$
we get $h_3(x)<0$ for any $x\geq 4$. Then by (\ref{r=3-e}) and (\ref{r=3-f}), we get
$$
g_3'(x)<0,\ \forall x\geq 4,
$$
which implies that  $g_3(x)$ is strictly decreasing on $[4,\infty)$. It follows that (\ref{r=3-c}) holds. The proof is complete.\hfill\fbox

%
%
%
%
%
%
%

\vskip 0.5cm

{ \noindent {\bf\large Acknowledgments}\quad We'd like to thank an anonymous referee for his/her useful suggestions to the first version of this manuscript.
 We thank Xuesong Li and Jingzi Yan for the discussing on the proof of Proposition 2.1 in the first version. This work was supported by the National Natural Science Foundation of China (12171335) and the Science Development Project of Sichuan University (2020SCUNL201).


\begin{thebibliography}{1234}

\bibitem{BPR23}  Bababesi, L.,  Pratelli, L.,  Rigo, P. 2023. On the Chv\'{a}tal-Janson conjecture. Statis. Probab. Lett.  194, 109744.

\bibitem{Do18} Doerr, B. 2018.  An elementary analysis of the probability that a binomial random variable exceeds its expectation. Statis. Probab. Lett.  139, 67-74.
\bibitem{GM14}  Greenberg, S.,  Mohri, M. 2014.  Tight lower bound on the probability of a binomial exceeding its expectation.   Statis. Probab. Lett. 86, 91-98.

\bibitem{Ja21}  Janson, J. 2021. On the probability that a binomial variable is at most its expectation.  Statis. Probab. Lett. 171, 109020.

\bibitem{PR16}  Pelekis, C.,  Ramon, J. 2016. A lower bound on the probability that a binomial random variable is exceeding its mean.  Statis. Probab. Lett. 119, 305-309.

\bibitem{Su21}  Sun, P.  Strictly unimodality of the probability that the binomial distribution is more than its expectation. Discrete Appl. Math.  301, 1-5.

\bibitem{SHS23}  Sun, P., Hu, Z.-C.,   Sun W. 2023. The extreme values of two probability functions for the Gamma distribution. arXiv: 2303.17487v1.

\bibitem{XLH22} Xu, K., Li., F.-B.,  Hu, Z.-C. 2022. Study on Poisson distribution and geometric distribution motivated by Chv\'{a}tal's conjecture. arXiv: 2210.16515v1.





\end{thebibliography}
\end{document}